\documentclass[english,3p]{elsarticle}
\usepackage[T1]{fontenc}
\usepackage[latin9]{inputenc}
\usepackage{amsmath}
\usepackage{amsthm}
\usepackage{amssymb}

\makeatletter
\theoremstyle{plain}
\newtheorem{thm}{\protect\theoremname}
  \theoremstyle{definition}
  \newtheorem{defn}[thm]{\protect\definitionname}
  \theoremstyle{plain}
  \newtheorem{lem}[thm]{\protect\lemmaname}
  \theoremstyle{plain}
  \newtheorem{prop}[thm]{\protect\propositionname}
  \theoremstyle{plain}
  \newtheorem{cor}[thm]{\protect\corollaryname}
  \theoremstyle{remark}
  \newtheorem{rem}[thm]{\protect\remarkname}

\makeatother

\usepackage{babel}
  \providecommand{\corollaryname}{Corollary}
  \providecommand{\definitionname}{Definition}
  \providecommand{\lemmaname}{Lemma}
  \providecommand{\propositionname}{Proposition}
  \providecommand{\remarkname}{Remark}
\providecommand{\theoremname}{Theorem}

\begin{document}

\begin{frontmatter}{}

\title{On the Well-Posedness of Two Driven-Damped Gross-Pitaevskii-Type
Models for Exciton-Polariton Condensates }

\author[rvt]{Jakob Möller\fnref{fn1}}

\author[rvt]{Jesus Sierra\corref{cor1}\fnref{fn2}}

\fntext[fn1]{jakob.moeller@univie.ac.at}

\fntext[fn2]{jesus.sierra@univie.ac.at }

\cortext[cor1]{Corresponding author}

\address[rvt]{Faculty of Mathematics, University of Vienna, Oskar-Morgenstern-Platz
1, 1090 Vienna, Austria.}
\begin{abstract}
We study the well-posedness of two systems modeling the non-equilibrium
dynamics of pumped decaying Bose-Einstein condensates. In particular,
we present the local theory for rough initial data using the Fourier
restricted norm method introduced by Bourgain. We extend the result
globally for initial data in $L^{2}$.
\end{abstract}

\end{frontmatter}{}

\section{Introduction}

In this paper, we study the (local) well-posedness theory of two closely
related models describing the (non-equilibrium) dynamics of pumped
decaying condensates, e.g., the Bose-Einstein condensation of exciton-polaritons.
The first model is the following driven-damped nonlinear Schrödinger
equation \cite{keeling2008spontaneous}:
\begin{equation}
i\partial_{t}u=-\partial_{x}^{2}u+\left|u\right|^{2}u+i\left(\xi-\sigma\left|u\right|^{2}\right)u,\label{eq:cGPE}
\end{equation}
where $u=u\left(x,t\right)$, $x\in\mathbb{T}$, $\xi,$ $\sigma$
are positive constants, and $u_{0}=u\left(x,0\right)\in H^{s}\left(\mathbb{T}\right)$,
$s\geq0$. 

The second model consists of a generalized open-dissipative Gross-Pitaevskii
equation for the macroscopic wave-function of the polaritons, $u=u\left(x,t\right)$,
coupled to a simple rate equation for the exciton reservoir density,
$n=n\left(x,t\right)$ \cite{wouters2007excitations,wouters2008spatial}:
\begin{align}
i\partial_{t}u= & -\partial_{x}^{2}u+g\left|u\right|^{2}u+\lambda nu+i\left(Rn-\alpha\right)u,\label{eq:SYS}\\
\partial_{t}n= & P-\left(R\left|u\right|^{2}+\beta\right)n,\nonumber 
\end{align}
subject to the initial data $\left.u\right|_{t=0}=u_{0}\left(x\right),$
$\left.n\right|_{t=0}=n_{0}\left(x\right)$, $x\in\mathbb{R}$. Above,
$\alpha,$ $\beta,$ $\lambda,$ $g,$ $R$ are positive constants
and $P=P\left(x\right)\geq0$ (compactly supported, bounded).

In our analysis, we shall consider (\ref{eq:cGPE}) on the one-dimensional
torus. This choice is physically motivated by the fact that a stable
condensate can only form in a spatially confined system. Such confinement
gives rise to some technical challenges due to the loss of dispersion.
Our approach is base on the Fourier restricted norm method introduced
by Bourgain in \cite{bourgain1993fourier,bourgain1993fourier2}. In
the case of the system (\ref{eq:SYS}) the confinement is given by
$P$. Our study of (\ref{eq:SYS}) requires some refinements of Bourgain's
method, in particular, the ones introduced by Kenig-Ponce-Vega in
\cite{kenig1996bilinear,kenig1996quadratic} and later used by Ginibre
et. al. in \cite{ginibre1997cauchy} to study the well-posedness theory
of the Zakharov system. On the other hand, it is important to notice
that (\ref{eq:SYS}) does not have derivatives in the nonlinearities.

We shall refer to (\ref{eq:cGPE}) as the complex Gross-Pitaevskii
equation and to (\ref{eq:SYS}) as the exciton-polariton system. 

\section{Well-posedness of the complex Gross-Pitaevskii equation}

Using Duhamel's principle, we consider the following integral equation
associated with (\ref{eq:cGPE}):
\begin{equation}
u\left(t\right)=S\left(t\right)u_{0}+\int_{0}^{t}S\left(t-\tau\right)\left(\xi u-\left(\sigma+i\right)\left|u\right|^{2}u\right)\left(\tau\right)d\tau,\label{eq:Duhamel_cGPE}
\end{equation}
where $S\left(t\right)=e^{it\partial_{x}^{2}}$. We introduce now
the basic notation and ideas related to the restricted norm method;
see, e.g., \cite{erdougan2016dispersive,ginibre1995probleme,tao2006nonlinear}
for a detailed review of this topic. 

Denote by $l_{k}^{q}L_{\tau}^{p}$ the Banach space $l_{k}^{q}\left(\mathbb{Z}:L_{\tau}^{p}\left(\mathbb{R}\right)\right)$.
Let $\hat{\cdot}$ stand for the Fourier transform with respect to
space-time, i.e.,
\[
\hat{g}\left(k,\tau\right)=\int_{-\infty}^{\infty}\int_{\mathbb{T}}\exp\left(-ikx-it\tau\right)g\left(x,t\right)dxdt.
\]
We denote by $\mathcal{F}_{x}$ the Fourier transform with respect
to the space variable
\[
\mathcal{F}_{x}g\left(k\right)=\int_{\mathbb{T}}\exp\left(-ikx\right)g\left(x\right)dx.
\]
\begin{defn}
Let $\mathcal{W}$ be the space of functions $u:\mathbb{T}\times\mathbb{R}\rightarrow\mathbb{C}$,
such that $u\left(x,\cdot\right)\in\mathcal{S}\left(\mathbb{R}\right)$
for each $x\in\mathbb{T}$ and $u\left(\cdot,t\right)\in C^{\infty}\left(\mathbb{T}\right)$
for each $t\in\mathbb{R}$. We define the space $X^{s,b}$ as the
completion of $\mathcal{W}$ with respect to the norm
\[
\left\Vert u\right\Vert _{X^{s,b}}=\left\Vert \left\langle k\right\rangle ^{s}\left\langle \tau+k^{2}\right\rangle ^{b}\hat{u}\left(k,\tau\right)\right\Vert _{l_{k}^{2}L_{\tau}^{2}}=\left\Vert e^{-it\partial_{x}^{2}}u\right\Vert _{H_{x}^{s}H_{t}^{b}},
\]
where $\left\langle \cdot\right\rangle =\left(1+\left|\cdot\right|^{2}\right)^{1/2}$
(Japanese bracket). 
\end{defn}

One can verify that the dual space of $X^{s,b}$ is $X^{-s,-b}$.
Moreover, $X^{s',b'}\subset X^{s,b}$ for $s'\geq s,$ $b'\geq b$.
Since we shall study the local theory using a contraction argument
in a time interval $\left[-\delta,\delta\right]$ with $\delta\leq1$,
we define the (restricted) space $X_{\delta}^{s,b}$ to be the equivalent
classes of functions that agree on $\left[-\delta,\delta\right]$,
with the norm 
\[
\left\Vert u\right\Vert _{X_{\delta}^{s,b}}=\underset{\tilde{u}=u,t\in\left[-\delta,\delta\right]}{\inf}\left\Vert \tilde{u}\right\Vert _{X^{s,b}}.
\]

Let $\eta\in C_{0}^{\infty}\left(\mathbb{R}\right)$ such that $\eta\left(t\right)=1$
for $t\in\left[-1,1\right]$. Define the operator 
\begin{equation}
\Gamma_{u_{0}}\left(u\right)=\eta\left(t\right)S\left(t\right)u_{0}+\eta\left(t\right)\int_{0}^{t}S\left(t-\tau\right)\left(\xi u-\left(\sigma+i\right)\left|u\right|^{2}u\right)\left(\tau\right)d\tau,\label{eq:fp_op_cGPE}
\end{equation}
on the ball
\begin{equation}
B_{R}=\left\{ u\in X_{\delta}^{s,b}:\left\Vert u\right\Vert _{X_{\delta}^{s,b}}\leq R\right\} ,\label{eq:Ball_cGPE}
\end{equation}
where $R=C\left\Vert u_{0}\right\Vert _{H^{s}}$, $s\geq0$. Note
that, since $\delta\leq1$, a fixed point of (\ref{eq:fp_op_cGPE})
gives a solution of the complex Gross-Pitaevskii equation on $\left[-\delta,\delta\right]$.
On the other hand, we have a constraint on the value of $b$ to ensure
continuity (in time) of these solutions, as the following lemma shows
(see \cite[Lemma 3.9]{erdougan2016dispersive}):
\begin{lem}
\label{lem:Emb_cont}For any $b>\frac{1}{2}$, $X_{\delta}^{s,b}\subset C_{t}^{0}H_{x}^{s}\left(\left[-\delta,\delta\right]\times\mathbb{T}\right)$.
\end{lem}

To handle our contraction argument, we shall use the following (see
\cite[Section 3.5.1]{erdougan2016dispersive}).
\begin{lem}
Let $0<\delta\leq1$, $s,b\in\mathbb{R}$. Then
\begin{equation}
\left\Vert \eta\left(t\right)e^{it\partial_{x}^{2}}u_{0}\right\Vert _{X_{\delta}^{s,b}}\leq C\left\Vert u_{0}\right\Vert _{H^{s}}.\label{eq:lemma3_a}
\end{equation}

For any $-\frac{1}{2}<b'<b<\frac{1}{2}$ and $s\in\mathbb{R},$ we
have
\begin{equation}
\left\Vert u\right\Vert _{X_{\delta}^{s,b'}}\leq C\delta^{b-b'}\left\Vert u\right\Vert _{X_{\delta}^{s,b}}.\label{Lemma3_b}
\end{equation}

Let $-\frac{1}{2}<b'\leq0$ and $b=b'+1$. Then
\begin{equation}
\left\Vert \eta\left(t\right)\int_{0}^{t}e^{i\left(t-s\right)\partial_{x}^{2}}F\left(s\right)ds\right\Vert _{X_{\delta}^{s,b}}\leq C\left\Vert F\right\Vert _{X_{\delta}^{s,b'}}.\label{eq:Lemma3_c}
\end{equation}
\end{lem}

The following result by Bourgain is essential for our analysis.
\begin{lem}
Let $u$ be a smooth space-time function. Then
\[
\left\Vert u\right\Vert _{L_{x\in\mathbb{T},t\in\mathbb{R}}^{4}}\leq C\left\Vert u\right\Vert _{X^{0,3/8}}.
\]
\end{lem}

Using the previous lemma, one can show the following (see \cite[Proposition 3.26]{erdougan2016dispersive}).
\begin{lem}
Let $s\geq0$. Then
\[
\left\Vert \left|u\right|^{2}u\right\Vert _{X_{\delta}^{s,-\frac{3}{8}}}\leq C\left\Vert u\right\Vert _{X_{\delta}^{0,\frac{3}{8}}}^{2}\left\Vert u\right\Vert _{X_{\delta}^{s,\frac{3}{8}}}.
\]
\end{lem}

Now we can present the main result of this section.
\begin{prop}
The complex Gross-Pitaevskii equation (\ref{eq:cGPE}) is locally
well-posed in $H_{x}^{s}\left(\mathbb{T}\right)$, $s\geq0$, i.e.,
for any $u_{0}\in H_{x}^{s}\left(\mathbb{T}\right)$ there is a unique
solution $u\in C_{t}^{0}H_{x}^{s}\left(\left[-\delta,\delta\right]\times\mathbb{T}\right)\cap X_{\delta}^{s,b}$,
with $\frac{1}{2}<b<\frac{5}{8}$. Moreover, the solution depends
continuously on the data.
\end{prop}

\begin{proof}
We run the contraction argument in $B_{R}\subset X_{\delta}^{s,b}$
(with $\frac{1}{2}<b<\frac{5}{8}$ and $\delta$ small enough) for
the operator $\Gamma_{u_{0}}\left(u\right)$ defined in (\ref{eq:fp_op_cGPE})-(\ref{eq:Ball_cGPE}).
Using (\ref{eq:lemma3_a}), (\ref{Lemma3_b}), (\ref{eq:Lemma3_c}),
and the embedding $X^{s',b'}\subset X^{s,b}$ for $s'\geq s,$ $b'\geq b$,
we obtain
\begin{align*}
\left\Vert \Gamma_{u_{0}}u\right\Vert _{X_{\delta}^{s,b}}\leq & C\left\Vert u_{0}\right\Vert _{H^{s}\left(\mathbb{T}\right)}+C\delta^{1-b-\frac{3}{8}}\left\Vert \xi u-\left(\sigma+i\right)\left|u\right|^{2}u\right\Vert _{X_{\delta}^{s,-\frac{3}{8}}}\\
\leq & C\left\Vert u_{0}\right\Vert _{H^{s}\left(\mathbb{T}\right)}+C\delta^{1-b-\frac{3}{8}}\left(\left\Vert u\right\Vert _{X_{\delta}^{s,-\frac{3}{8}}}+\left\Vert \left|u\right|^{2}u\right\Vert _{X_{\delta}^{s,-\frac{3}{8}}}\right)\\
\leq & C\left\Vert u_{0}\right\Vert _{H^{s}\left(\mathbb{T}\right)}+C\delta^{1-b-\frac{3}{8}}\left(\left\Vert u\right\Vert _{X_{\delta}^{s,-\frac{3}{8}}}+\left\Vert u\right\Vert _{X_{\delta}^{0,\frac{3}{8}}}^{2}\left\Vert u\right\Vert _{X_{\delta}^{s,\frac{3}{8}}}\right)\\
\leq & C\left\Vert u_{0}\right\Vert _{H^{s}\left(\mathbb{T}\right)}+C\delta^{1-b-\frac{3}{8}}\left(\left\Vert u\right\Vert _{X_{\delta}^{s,b}}+\left\Vert u\right\Vert _{X_{\delta}^{s,b}}^{3}\right).
\end{align*}
Similar estimates hold for the difference. We omit the standard details.
Note that, since $b>\frac{1}{2}$, by Lemma \ref{lem:Emb_cont} the
solution is continuous in time with values in $H^{s}\left(\mathbb{T}\right)$,
$s\geq0$.
\end{proof}
\begin{cor}
The complex Gross-Pitaevskii equation (\ref{eq:cGPE}) is globally
well-posed in $L_{x}^{2}\left(\mathbb{T}\right)$.
\end{cor}

\begin{proof}
Multiply (\ref{eq:cGPE}) by $\bar{u}$, take the imaginary part,
and use integration by parts to obtain
\begin{equation}
\frac{d}{dt}\int_{\mathbb{T}}\left|u\right|^{2}dx-2\xi\int_{\mathbb{T}}\left|u\right|^{2}dx+2\sigma\int\left|u\right|^{4}dx=0.\label{eq:F1}
\end{equation}
Since $\left(\sqrt{\sigma}s^{2}-\frac{2\xi}{\sqrt{\sigma}}\right)^{2}\geq0$,
we have $\sigma s^{4}-4\xi s^{2}+\frac{4\xi^{2}}{\sigma}\geq0$, for
$s\in\mathbb{R}.$ Setting $s=\left|u\right|$ and integrating over
$\mathbb{T}$ yield
\begin{equation}
-\sigma\int_{\mathbb{T}}\left|u\right|^{4}dx+4\xi\int_{\mathbb{T}}\left|u\right|^{2}dx\leq\frac{4\xi^{2}}{\sigma}\left|\mathbb{T}\right|,\label{eq:F2}
\end{equation}
where $\left|\mathbb{T}\right|$ is the measure of $\mathbb{T}$.
Combining (\ref{eq:F1}) and (\ref{eq:F2}) gives
\[
\frac{d}{dt}\int_{\mathbb{T}}\left|u\right|^{2}dx+2\xi\int_{\mathbb{T}}\left|u\right|^{2}dx+\sigma\int\left|u\right|^{4}dx\leq\frac{4\xi^{2}}{\sigma}\left|\mathbb{T}\right|.
\]
Using the last expression along with Gronwall's lemma, we obtain
\begin{equation}
\left\Vert u\left(t\right)\right\Vert _{L^{2}}^{2}\leq\left\Vert u_{0}\right\Vert _{L^{2}}^{2}e^{-2\xi t}+\frac{2\xi}{\sigma}\left|\mathbb{T}\right|\left(1-e^{-2\xi t}\right),\textrm{ }t\geq0.\label{eq:Abs_set}
\end{equation}
From (\ref{eq:Abs_set}) we conclude that the local solution in $L_{x}^{2}\left(\mathbb{T}\right)$
can be extended globally. Note that to justify the calculations above
we need to use continuous dependence on the data, approximate $u_{0}$
by a sequence of smooth functions, and take the limit. See, e.g.,
\cite{cazenave2003semilinear,linares2014introduction} for a detailed
description of this procedure.
\end{proof}
\begin{rem}
Letting $t\rightarrow\infty$ in (\ref{eq:Abs_set}) gives
\[
\underset{t\rightarrow\infty}{\limsup}\left\Vert u\right\Vert _{L^{2}}^{2}\leq\frac{2\xi}{\sigma}\left|\mathbb{T}\right|,
\]
which guarantees the existence of an absorbing set for the complex
Gross-Pitaevskii equation in $L^{2}\left(\mathbb{T}\right)$. See,
e.g., \cite{temam2012infinite}. 

The $H^{1}$ theory and stationary solutions of the complex Gross-Pitaevskii
equation (in the full domain with a harmonic trapping potential and
$\xi$ space dependent with compact support) have been studied in
\cite{feng2018global,hajaiej2021ground,sierra2015gross}.
\end{rem}

\section{Well-posedness of the exciton-polariton system }

Using Duhamel's principle, we consider the following integral equations
associated with (\ref{eq:SYS}):
\begin{equation}
u\left(t\right)=S\left(t\right)u_{0}+\int_{0}^{t}S\left(t-s\right)\left(-ig\left|u\right|^{2}u+\left(R-i\lambda\right)nu-\alpha u\right)\left(s\right)ds,\label{eq:EP1}
\end{equation}
\begin{equation}
n\left(t\right)=n_{0}+\int_{0}^{t}\left(P-R\left|u\right|^{2}n-\beta n\right)\left(s\right)ds,\label{eq:EP2}
\end{equation}
where $S\left(t\right)=e^{it\partial_{x}^{2}}$. Our choice for the
second expression is because the corresponding equation in (\ref{eq:SYS})
is an ODE in $n$; hence, it does not have an appropriate dispersion
relation for the subsequent analysis. 
\begin{defn}
Let $X^{s,b}$ be the Banach space of functions on $\mathbb{R}\times\mathbb{R}$
defined by the norm
\[
\left\Vert u\right\Vert _{X_{\phi\left(\xi\right)}^{s,b}}=\left\Vert \left\langle \xi\right\rangle ^{s}\left\langle \tau+\phi\left(\xi\right)\right\rangle ^{b}\hat{u}\left(\xi,\tau\right)\right\Vert _{L_{\xi.\tau}^{2}},
\]
where $\phi$ corresponds to the dispersion relation of the equation
under consideration. We usually write $\left\Vert \cdot\right\Vert _{2}=\left\Vert \cdot\right\Vert _{L_{\xi,\tau}^{2}}$.
\end{defn}

Similarly, we define the auxiliary spaces $Y^{s}$ by the norm
\[
\left\Vert u\right\Vert _{Y_{\phi\left(\xi\right)}^{s}}=\left\Vert \left\langle \xi\right\rangle ^{s}\left\langle \tau+\phi\left(\xi\right)\right\rangle ^{-1}\hat{u}\left(\xi,\tau\right)\right\Vert _{L_{\xi}^{2}L_{\tau}^{1}}.
\]

We want to solve the Cauchy problem corresponding to (\ref{eq:SYS})
in the context of the previous spaces and in some time interval $\left[-T,T\right]$.
To reach this goal, it is convenient to introduce a cutoff in (\ref{eq:EP1})-(\ref{eq:EP2}).
Let $\psi\in C_{0}^{\infty}\left(\mathbb{R}\right)$ be even, $0\leq\psi\leq1$,
such that $\psi=1$ on $\left[-1,1\right]$ and $\textrm{supp }\psi\subset\left(-2,2\right)$.
Furthermore, let $\psi_{T}\left(t\right)=\psi\left(t/T\right)$, $0<T\leq1$.
The cutoff version of (\ref{eq:EP1})-(\ref{eq:EP2}) is given by
\begin{equation}
u\left(t\right)=\psi\left(t\right)S\left(t\right)u_{0}+\psi_{T}\left(t\right)\int_{0}^{t}S\left(t-s\right)\left(-ig\left|u\right|^{2}u+\left(R-i\lambda\right)nu-\alpha u\right)\left(s\right)ds,\label{eq:EP1_cut}
\end{equation}
\begin{equation}
n\left(t\right)=\psi\left(t\right)n_{0}+\psi_{T}\left(t\right)\int_{0}^{t}\left(P-R\left|u\right|^{2}n-\beta n\right)\left(s\right)ds.\label{eq:EP2_cut}
\end{equation}

Like in the previous section, we define the restricted space $\left\Vert u\right\Vert _{X_{\phi\left(\xi\right),T}^{s,b}}$
as the equivalent classes of functions that agree on $t\in\left[-T,T\right]$,
with the norm
\[
\left\Vert u\right\Vert _{X_{\phi\left(\xi\right),T}^{s,b}}=\underset{\tilde{u}=u,t\in\left[-T,T\right]}{\inf}\left\Vert \tilde{u}\right\Vert _{X_{\phi\left(\xi\right)}^{s,b}}.
\]
Similarly, we define the space $\left\Vert u\right\Vert _{Y_{\phi\left(\xi\right),T}^{s}}$.

The following lemma will be the starting point for our contraction
argument (see \cite[Lemma 2.1]{ginibre1997cauchy}).
\begin{lem}
\label{lem:Lem_ini}Let $s\in\mathbb{R}$, $b'\leq0\leq b\leq b'+1$,
and $T\leq1$. Then
\[
\left\Vert \psi_{T}\left(t\right)\int_{0}^{t}S\left(t-\tau\right)F\left(\tau\right)d\tau\right\Vert _{X_{\phi\left(\xi\right)=\xi^{2}}^{s,b}}\leq C\left(T^{1-b+b'}\left\Vert F\right\Vert _{X_{\phi\left(\xi\right)=\xi^{2}}^{s,b'}}+T^{1/2-b}\left\Vert F\right\Vert _{Y_{\phi\left(\xi\right)=\xi^{2}}^{s}}\right),\textrm{ }
\]
\[
\left\Vert \psi_{T}\left(t\right)\int_{0}^{t}F\left(\tau\right)d\tau\right\Vert _{X_{\phi\left(\xi\right)\equiv0}^{s,b}}\leq C\left(T^{1-b+b'}\left\Vert F\right\Vert _{X_{\phi\left(\xi\right)\equiv0}^{s,b'}}+T^{1/2-b}\left\Vert F\right\Vert _{Y_{\phi\left(\xi\right)\equiv0}^{s}}\right).
\]
Furthermore, if $b'>-1/2$,
\[
\left\Vert \psi_{T}\left(t\right)\int_{0}^{t}S\left(t-\tau\right)F\left(\tau\right)d\tau\right\Vert _{X_{\phi\left(\xi\right)=\xi^{2}}^{s,b}}\leq CT^{1-b+b'}\left\Vert F\right\Vert _{X_{\phi\left(\xi\right)=\xi^{2}}^{s,b'}},
\]
\[
\left\Vert \psi_{T}\left(t\right)\int_{0}^{t}F\left(\tau\right)d\tau\right\Vert _{X_{\phi\left(\xi\right)\equiv0}^{s,b}}\leq CT^{1-b+b'}\left\Vert F\right\Vert _{X_{\phi\left(\xi\right)\equiv0}^{s,b'}}.
\]
\end{lem}

As mentioned before, $X_{\phi\left(\xi\right)}^{s,b}\subset C\left(\mathbb{R},H^{s}\right)$,
$b>1/2$. This is no longer valid if $b\leq1/2$, and this is why
we need to consider the spaces $Y_{\phi\left(\xi\right)}^{s}$ (see
\cite[Lemma 2.2]{ginibre1997cauchy}).

We now follow closely the ideas presented in \cite{ginibre1997cauchy}.
As mentioned before, we want to solve the cutoff integral version
of the exciton-polariton system (\ref{eq:EP1_cut})-(\ref{eq:EP2_cut})
by a contraction method with $u\in X_{\phi\left(\xi\right)=\xi^{2}}^{k,a_{2}}$
and $n\in X_{\phi\left(\xi\right)\equiv0}^{l,a}$ for suitable $a$,
$a_{2}$, and $k$, $l$. We start by estimating the nonlinearity
\[
f_{1}=nu
\]
in $X_{\phi\left(\xi\right)=\xi^{2}}^{k,-a_{1}}$ for suitable $a_{1}$.

We estimate $\hat{f}_{1}\left(\xi_{1},\tau_{1}\right)$ in terms of
$\hat{n}\left(\xi,\tau\right)$ and $\hat{u}\left(\xi_{2},\tau_{2}\right)$.
We have the following relations due to the convolution structure 
\[
\xi=\xi_{1}-\xi_{2},
\]
\[
\tau=\tau_{1}-\tau_{2}.
\]
We also introduce the variables
\[
\sigma_{1}=\tau_{1}+\xi_{1}^{2},
\]
\[
\sigma_{2}=\tau_{2}+\xi_{2}^{2},
\]
\[
\sigma=\tau.
\]
In terms of these variables, we have
\begin{equation}
z\equiv\xi_{1}^{2}-\xi_{2}^{2}=\sigma_{1}-\sigma_{2}-\sigma.\label{eq:Def_z}
\end{equation}
We use this expression to obtain estimates of $\xi_{1}^{2}$ (resp.
$\xi_{2}^{2}$) in terms of $\xi_{2}^{2}$ (resp. $\xi_{1}^{2}$)
and of the $\sigma\textrm{'s}$.

To estimate $f_{1}$, we define $\hat{v}_{2}=\left\langle \xi_{2}\right\rangle ^{k}\left\langle \sigma_{2}\right\rangle ^{a_{2}}\hat{u}$
and $\hat{v}=\left\langle \xi\right\rangle ^{l}\left\langle \sigma\right\rangle ^{a}\hat{n}$
so that
\[
\left\Vert u\right\Vert _{X_{\phi\left(\xi\right)=\xi^{2}}^{k,a_{2}}}=\left\Vert v_{2}\right\Vert _{2},
\]
and
\[
\left\Vert n\right\Vert _{X_{\phi\left(\xi\right)\equiv0}^{l,a}}=\left\Vert v\right\Vert _{2}.
\]
To estimate $f_{1}$ in $X_{\phi\left(\xi\right)=\xi^{2}}^{k,-a_{1}}$,
we take its scalar product with a generic function in $X_{\phi\left(\xi\right)=\xi^{2}}^{-k,a_{1}}$
with Fourier transform $\left\langle \xi_{1}\right\rangle ^{k}\left\langle \sigma_{1}\right\rangle ^{-a_{1}}\hat{v}_{1}$
and $v_{1}\in L^{2}$. Then the required estimate in $X_{\phi\left(\xi\right)=\xi^{2}}^{k,-a_{1}}$
takes the form
\begin{equation}
\left|S\right|\leq C\left\Vert v\right\Vert _{2}\left\Vert v_{1}\right\Vert _{2}\left\Vert v_{2}\right\Vert _{2},\label{eq:main_est}
\end{equation}
where
\begin{equation}
S=\int\frac{\hat{v}\hat{v}_{1}\hat{v}_{2}\left\langle \xi_{1}\right\rangle ^{k}}{\left\langle \sigma\right\rangle ^{a}\left\langle \sigma_{1}\right\rangle ^{a_{1}}\left\langle \sigma_{2}\right\rangle ^{a_{2}}\left\langle \xi_{2}\right\rangle ^{k}\left\langle \xi\right\rangle ^{l}},\label{eq:def_S}
\end{equation}
and
\[
\hat{v}=\hat{v}\left(\xi,\tau\right),
\]
\[
\hat{v}_{1}=\hat{v}_{1}\left(\xi_{1},\tau_{1}\right),
\]
\[
\hat{v}_{2}=\hat{v}_{2}\left(\xi_{2},\tau_{2}\right),
\]
constrained by
\[
\xi=\xi_{1}-\xi_{2},\textrm{ }\tau=\tau_{1}-\tau_{2},
\]
and the integral is over $d\xi_{1}$, $d\xi_{2}$, $d\tau_{1}$, $d\tau_{2}$. 

We often use the following two elementary facts in our analysis.
\begin{lem}
\label{lem:Conv}Let $f\in L^{q}\left(\mathbb{R}\right)$, $g\in L^{q'}\left(\mathbb{R}\right),$
$1\leq q,q'\leq\infty,$ $1/q+1/q'=1.$ Assume that $f$ and $g$
are nonnegative, even, and non-increasing for positive argument. Then,
$f*g$ has the same properties.
\end{lem}

One can use Lemma \ref{lem:Conv} to show that $f*g$ takes its maximum
at zero. Using this fact, we can show the following
\begin{lem}
\label{lem:bracket_lemma}Let $0\leq a_{-}\leq a_{+}$ and $a_{+}+a_{-}>1/2,$
then the following estimate holds for all $s\in\mathbb{R}$
\[
J\left(s\right)=\int\left\langle y-s\right\rangle ^{-2a_{+}}\left\langle y+s\right\rangle ^{-2a_{-}}dy\leq C\left\langle s\right\rangle ^{-\alpha},
\]
where $\alpha=2a_{-}-\left[1-2a_{+}\right]_{+}$.
\end{lem}

See \cite{ginibre1997cauchy} for a proof of the previous lemmata.
\begin{lem}
\label{lem:Lemm_S}Let $k,$ $l,$ $a,$ $a_{1},$ $a_{2}$ satisfy
\begin{equation}
\begin{array}{ccc}
l\geq-1/2, & k\geq0, & k-l\leq1,\end{array}\label{eq:Lem1_a}
\end{equation}
\begin{equation}
\begin{array}{ccc}
a,a_{1},a_{2}>1/4, & a+a_{1}>3/4, & a+a_{2}>3/4,\end{array}\label{eq:Lem1_b}
\end{equation}
\begin{equation}
k-l\leq2a_{1},\label{eq:Lem1_c}
\end{equation}
then the estimate (\ref{eq:main_est}) holds.
\end{lem}

\begin{proof}
The principle of the proof is the following application of the Schwarz
inequality. Let $\zeta=\left(\xi,\tau\right),$ $\zeta_{i}=\left(\xi_{i},\tau_{i}\right),$
$i=1,2$ so that $\zeta=\zeta_{1}-\zeta_{2}.$ We want to estimate
an integral of the form
\[
J=\int\hat{v}\left(\zeta\right)\hat{v}_{1}\left(\zeta_{1}\right)\hat{v}_{2}\left(\zeta_{2}\right)K\left(\zeta_{1},\zeta_{2}\right)d\zeta_{1}d\zeta_{2}.
\]
Note that $\zeta_{1}=\zeta+\zeta_{2}$. Then, considering the Schwarz
inequality with respect to $\zeta$ we obtain
\begin{align*}
\left|J\right|^{2}\leq & \left\Vert v\right\Vert _{2}^{2}\int d\zeta\left|\int\hat{v}_{1}\left(\zeta+\zeta_{2}\right)\hat{v}_{2}\left(\zeta_{2}\right)K\left(\zeta+\zeta_{2},\zeta_{2}\right)d\zeta_{2}\right|^{2}\\
 & \textrm{(Schwarz w.r.t. }\zeta_{2}\textrm{ and extract sup)}\\
\leq & \left\Vert v\right\Vert _{2}^{2}\left\{ \underset{\zeta}{\sup}\int\left|K\left(\zeta+\zeta_{2},\zeta_{2}\right)\right|^{2}d\zeta_{2}\right\} \int\left|\hat{v}_{1}\left(\zeta+\zeta_{2}\right)\hat{v}_{2}\left(\zeta_{2}\right)\right|^{2}d\zeta d\zeta_{2}\\
 & \textrm{(use Fubini, translation invariance, and Plancherel)}\\
= & C^{2}\left\Vert v\right\Vert _{2}^{2}\left\Vert v_{1}\right\Vert _{2}^{2}\left\Vert v_{2}\right\Vert _{2}^{2},
\end{align*}
with
\begin{equation}
C^{2}=\underset{\zeta}{\sup}\int_{\zeta}\left|K\left(\zeta_{1},\zeta_{2}\right)\right|^{2}d\zeta_{2},\label{eq:Def_C}
\end{equation}
and the last integral runs over $\zeta_{2}$ (or $\zeta_{1}$) for
fixed $\zeta$. One obtains two similar estimates by circularly permuting
the variables and functions 1, 2, and 1-2 (the ones with no subindex).

Moreover, we define
\[
\alpha=2\min\left(a_{1},a_{2}\right)-\left[1-2\max\left(a_{1},a_{2}\right)\right]_{+},
\]
\[
\alpha_{1}=2\min\left(a,a_{2}\right)-\left[1-2\max\left(a,a_{2}\right)\right]_{+},
\]
\[
\alpha_{2}=2\min\left(a,a_{1}\right)-\left[1-2\max\left(a,a_{1}\right)\right]_{+}.
\]
We start by considering a particular case for $k$ and $l$. 

\[
\mathbf{Case}\textrm{ }k=0,\textrm{ }l=-1/2.
\]

In this case, the factors containing the $\xi$'s reduce to $\left\langle \xi\right\rangle ^{1/2}.$
Note that
\[
\left\langle \xi\right\rangle \leq1+\left|\xi\right|,
\]
then
\[
\left\langle \xi\right\rangle ^{1/2}\leq\left(1+\left|\xi\right|\right)^{1/2}\leq1+\left|\xi\right|^{1/2}.
\]
Therefore, for this case
\begin{equation}
S\leq\int\frac{\left|\hat{v}\hat{v}_{1}\hat{v}_{2}\right|}{\left\langle \sigma\right\rangle ^{a}\left\langle \sigma_{1}\right\rangle ^{a_{1}}\left\langle \sigma_{2}\right\rangle ^{a_{2}}}+\int\frac{\left|\hat{v}\hat{v}_{1}\hat{v}_{2}\right|\left|\xi\right|^{1/2}}{\left\langle \sigma\right\rangle ^{a}\left\langle \sigma_{1}\right\rangle ^{a_{1}}\left\langle \sigma_{2}\right\rangle ^{a_{2}}}=:A+Z_{0}.\label{eq:s_leq_A_Z0}
\end{equation}
Lemma \ref{lem:Lemm_ver_A} gives the bounds for $A$. For $Z_{0}$
we consider the following subregions.

$Region\textrm{ }\sigma\textrm{ }dominant,\textrm{ }i.e.,\textrm{ }\left|\sigma\right|\geq\max\left(\left|\sigma_{1}\right|,\left|\sigma_{2}\right|\right).$
We use directly (\ref{eq:Def_C}) and obtain
\[
C_{*}^{2}=\underset{\xi,\sigma}{\sup}\left\langle \sigma\right\rangle ^{-2a}\int_{*}\left|\xi\right|\left\langle \sigma_{1}\right\rangle ^{-2a_{1}}\left\langle \sigma_{2}\right\rangle ^{-2a_{2}}d\xi_{2}d\sigma_{2},
\]
where the integral is taken at fixed $\xi,$ $\sigma$. Now for fixed
$\xi,$ $\sigma,$ and $\sigma_{2}$, it follows from (\ref{eq:Def_z})
that
\[
2\left|\xi\right|d\xi_{2}=dz=d\sigma_{1},
\]
since $z=\xi_{1}^{2}-\xi_{2}^{2}$ and $\xi=\xi_{1}-\xi_{2}\Rightarrow\xi_{1}=\xi+\xi_{2}$,
which gives $dz=2\xi d\xi_{2}$. Therefore,
\[
C_{*}^{2}\leq C\underset{\sigma}{\sup}\left\langle \sigma\right\rangle ^{-2a}\int_{0}^{\left|\sigma\right|}\left\langle \sigma_{1}\right\rangle ^{-2a_{1}}d\sigma_{1}\int_{0}^{\left|\sigma\right|}\left\langle \sigma_{2}\right\rangle ^{-2a_{2}}d\sigma_{2}.
\]
Note that $\left\langle u\right\rangle \geq\left|u\right|,$ hence
$\left\langle u\right\rangle ^{-2a_{1}}\leq\left|u\right|^{-2a_{1}}.$
Then,
\begin{align*}
\int_{0}^{\left|\sigma\right|}\left\langle \sigma_{1}\right\rangle ^{-2a_{1}}d\sigma_{1}\leq & \int_{0}^{\left|\sigma\right|}\sigma_{1}^{-2a_{1}}d\sigma_{1}\\
\leq & C\left|\sigma\right|^{\left[1-2a_{1}\right]_{+}}\\
\leq & C\left\langle \sigma\right\rangle ^{\left[1-2a_{1}\right]_{+}}.
\end{align*}
Hence
\[
C_{*}^{2}\leq C\underset{\sigma}{\sup}\left\langle \sigma\right\rangle ^{-2a+\left[1-2a_{1}\right]_{+}+\left[1-2a_{2}\right]_{+}}.
\]
The last quantity is finite provided
\[
2a-\left[1-2a_{1}\right]_{+}-\left[1-2a_{2}\right]_{+}\geq0,
\]
which holds under the conditions
\[
a>0,\textrm{ }a_{1}+a>1/2,\textrm{ }a_{2}+a>1/2,\textrm{ }a+a_{1}+a_{2}>1.
\]

$Region\textrm{ }\sigma_{1}\textrm{ }dominant,\textrm{ }i.e.,\textrm{ }\left|\sigma_{1}\right|\geq\max\left(\left|\sigma\right|,\left|\sigma_{2}\right|\right).$
We now use the analog of (\ref{eq:Def_C}) with fixed $\zeta_{1}$
and obtain
\begin{equation}
C_{1}^{2}=\underset{\xi_{1},\sigma_{1}}{\sup}\left\langle \sigma_{1}\right\rangle ^{-2a_{1}}\int_{1}\left|\xi\right|\left\langle \sigma\right\rangle ^{-2a}\left\langle \sigma_{2}\right\rangle ^{-2a_{2}}d\xi_{2}d\sigma_{2},\label{eq:C1_sp1}
\end{equation}
where the integral is taken at fixed $\xi_{1},$ $\sigma_{1}$. To
continue the estimate, we split the $\sigma_{1}$ dominant region
into two subregions.

$Subregion\textrm{ }\left|\xi_{1}\right|\leq2\left|\xi_{2}\right|.$
Recall that $\xi=\xi_{1}-\xi_{2}$, hence, $\left|\xi\right|=\left|\xi_{1}-\xi_{2}\right|\leq\left|\xi_{1}\right|+\left|\xi_{2}\right|\leq3\left|\xi_{2}\right|$,
the last inequality due to the subregion. Furthermore, for fixed $\xi_{1},$
$\sigma_{1},$ and $\sigma_{2},$ it follows from (\ref{eq:Def_z})
that $2\left|\xi_{2}\right|d\xi_{2}=dz=d\sigma.$ Therefore,
\begin{align*}
C_{1}^{2}\leq & C\underset{\sigma_{1}}{\sup}\left\langle \sigma_{1}\right\rangle ^{-2a_{1}}\int_{0}^{\left|\sigma_{1}\right|}\left\langle \sigma\right\rangle ^{-2a}d\sigma\int_{0}^{\left|\sigma_{1}\right|}\left\langle \sigma_{2}\right\rangle ^{-2a_{2}}d\sigma_{2}\\
\leq C & \underset{\sigma_{1}}{\sup}\left\langle \sigma_{1}\right\rangle ^{-2a_{1}+\left[1-2a\right]_{+}+\left[1-2a_{2}\right]_{+}}.
\end{align*}
The last inequality is finite provided
\[
2a_{1}-\left[1-2a\right]_{+}-\left[1-2a_{2}\right]_{+}\geq0,
\]
which holds when
\[
a_{1}>0,\textrm{ }a+a_{1}>1/2,\textrm{ }a_{1}+a_{2}>1/2,\textrm{ }a+a_{1}+a_{2}>1.
\]

$Subregion\textrm{ }\left|\xi_{1}\right|\geq2\left|\xi_{2}\right|.$
In this region, note that
\begin{align*}
\left|\xi_{1}\right|\geq & 2\left|\xi_{2}\right|=2\left|\xi_{1}-\xi\right|\geq2\left|\left|\xi_{1}\right|-\left|\xi\right|\right|\\
= & 2\left|\left|\xi\right|-\left|\xi_{1}\right|\right|\geq2\left(\left|\xi\right|-\left|\xi_{1}\right|\right).
\end{align*}
Then
\begin{equation}
3\left|\xi_{1}\right|\geq2\left|\xi\right|\Rightarrow\left|\xi\right|\leq\frac{3}{2}\left|\xi_{1}\right|.\label{eq:xi_sub_sig1}
\end{equation}
Moreover
\[
\left|\xi_{1}\right|\geq2\left|\xi_{2}\right|\Rightarrow\xi_{1}^{2}\geq4\xi_{2}^{2}\Rightarrow-\xi_{1}^{2}\leq-4\xi_{2}^{2}\Rightarrow-\frac{1}{4}\xi_{1}^{2}\leq-\xi_{2}^{2}.
\]
Combining the last expression with (\ref{eq:Def_z}) and the fact
that we are in the region $\sigma_{1}$ dominant, we obtain
\[
\frac{3}{4}\xi_{1}^{2}=\xi_{1}^{2}-\frac{1}{4}\xi_{1}^{2}\leq\xi_{1}^{2}-\xi_{2}^{2}=\sigma_{1}-\sigma_{2}-\sigma\leq3\left|\sigma_{1}\right|,
\]
and therefore
\begin{equation}
\xi_{1}^{2}\leq4\left|\sigma_{1}\right|.\label{eq:xi1_sub_sig1}
\end{equation}
By (\ref{eq:xi_sub_sig1}), $\left|\xi\right|\leq C\left\langle \xi_{1}\right\rangle $,
and by (\ref{eq:xi1_sub_sig1}), $\left\langle \sigma_{1}\right\rangle ^{-2a_{1}}\leq C\left\langle \xi_{1}\right\rangle ^{-4a_{1}}$.
Using these facts and taking $y=\xi_{2}^{2}$ as integration variable
instead of $\xi_{2}$, we obtain
\begin{equation}
C_{1}^{2}\leq\underset{\xi_{1},\sigma_{1}}{\sup}\left\langle \xi_{1}\right\rangle ^{1-4a_{1}}\int_{0}^{\xi_{1}^{2}/4}y^{-1/2}dy\int\left\langle \sigma\right\rangle ^{-2a}\left\langle \sigma_{2}\right\rangle ^{-2a_{2}}d\sigma_{2},\label{eq:C1_sp2}
\end{equation}
where the boundary of the first integral is due to the subregion that
we are considering: $\left|\xi_{1}\right|\geq2\left|\xi_{2}\right|\Rightarrow\xi_{2}^{2}\leq\frac{1}{4}\xi_{1}^{2}$.
Note that, since
\[
\xi_{1}^{2}-\xi_{2}^{2}=\sigma_{1}-\sigma_{2}-\sigma\textrm{ }(\xi_{2}^{2}=y),
\]
then
\[
\xi_{1}^{2}-y-\sigma_{1}=-\sigma_{2}-\sigma\Rightarrow\left\langle \sigma_{2}+\left(\xi_{1}^{2}-y-\sigma_{1}\right)\right\rangle =\left\langle -\sigma\right\rangle =\left\langle \sigma\right\rangle .
\]
Hence
\[
C_{1}^{2}\leq C\underset{\xi_{1},\sigma_{1}}{\sup}\left\langle \xi_{1}\right\rangle ^{1-4a_{1}}\int_{0}^{\xi_{1}^{2}}y^{-1/2}dy\int\left\langle \sigma_{2}+\left(\xi_{1}^{2}-y-\sigma_{1}\right)\right\rangle ^{-2a}\left\langle \sigma_{2}\right\rangle ^{-2a_{2}}d\sigma_{2}.
\]
We estimate the last integral for fixed $\xi_{1},$ $\sigma,$ $\xi_{2},$
by Lemma \ref{lem:bracket_lemma}. Then
\[
C_{1}^{2}\leq C\underset{\xi_{1},\sigma_{1}}{\sup}\left\langle \xi_{1}\right\rangle ^{1-4a_{1}}\int_{0}^{\xi_{1}^{2}}\left\langle \xi_{1}^{2}-y-\sigma_{1}\right\rangle ^{-\alpha_{1}}y^{-1/2}dy.
\]
We extend the range of integration of $y$ symmetrically to $\left[-\xi_{1}^{2}/4,\xi_{1}^{2}/4\right]$
and apply Lemma \ref{lem:Conv} with $f\left(y\right)=\left|y\right|^{-1/2}\chi\left(\left|y\right|\leq\xi_{1}^{2}/4\right)$
and $g\left(y\right)=\left\langle \xi_{1}^{2}-y-\sigma_{1}\right\rangle ^{-\alpha_{1}}$
to conclude that the supremum over $\sigma_{1}$ is attained for $\sigma_{1}=\xi_{1}^{2}$.
Hence,
\[
C_{1}^{2}\leq C\underset{\xi_{1}}{\sup}\left\langle \xi_{1}\right\rangle ^{1-4a_{1}}\int_{0}^{\xi_{1}^{2}}\left\langle y\right\rangle ^{-\alpha_{1}}y^{-1/2}dy.
\]
The last quantity is finite, provided $a_{1}\geq1/4$ and $\alpha_{1}>1/2.$
The latter is equivalent to
\[
a>1/4,\textrm{ }a_{2}>1/4,\textrm{ }a+a_{2}>3/4.
\]

$Region\textrm{ }\sigma_{2}\textrm{ }dominant.$ This region is obtained
from the previous one by exchanging 1 and 2. This has the effect of
exchanging $a_{2}$ and $a_{1},$ so that the same proof applies since
the only assumption used so far, namely (\ref{eq:Lem1_b}), is symmetric
in $a_{2}$ and $a_{1}$.

\[
\mathbf{General}\textrm{ }k\textrm{ }\mathbf{and}\textrm{ }l,\textrm{ }k\geq0
\]

We consider separately the regions $\left|\xi_{1}\right|\leq2\left|\xi_{2}\right|$
and $\left|\xi_{1}\right|\geq2\left|\xi_{2}\right|$.
\[
Region\textrm{ }\left|\xi_{1}\right|\leq2\left|\xi_{2}\right|
\]

In this region
\[
\left\langle \xi_{1}\right\rangle ^{k}\left\langle \xi_{2}\right\rangle ^{-k}\left\langle \xi\right\rangle ^{-l}\leq C\left\langle \xi\right\rangle ^{-l},
\]
so that the factors with $k$'s disappear and the resulting expression
is decreasing in $l$. It is therefore sufficient to derive estimate
(\ref{eq:main_est}) in the case $l=-1/2$, which is the special case
considered previously.
\[
Region\textrm{ }\left|\xi_{1}\right|\geq2\left|\xi_{2}\right|
\]

In this region, we have
\[
\left|\xi_{1}\right|\geq2\left|\xi_{2}\right|\Rightarrow-\left|\xi_{2}\right|\geq-\frac{1}{2}\left|\xi_{1}\right|,
\]
and hence
\[
\left|\xi\right|=\left|\xi_{1}-\xi_{2}\right|\geq\left|\left|\xi_{1}\right|-\left|\xi_{2}\right|\right|\geq\frac{1}{2}\left|\xi_{1}\right|\Rightarrow\left|\xi_{1}\right|\leq2\left|\xi\right|.
\]
Moreover, from (\ref{eq:xi_sub_sig1}) we have $3\left|\xi_{1}\right|\geq2\left|\xi\right|$.
Therefore,
\[
\left|\xi_{1}\right|\leq2\left|\xi\right|\leq3\left|\xi_{1}\right|.
\]
We deduce
\begin{equation}
\left\langle \xi_{1}\right\rangle \leq C_{1}\left\langle \xi\right\rangle \leq C_{2}\left\langle \xi_{1}\right\rangle .\label{eq:bound_xi_xi1}
\end{equation}
Now, using (\ref{eq:bound_xi_xi1}), we get
\begin{align*}
\int\frac{\hat{v}\hat{v}_{1}\hat{v}_{2}\left\langle \xi_{1}\right\rangle ^{k}}{\left\langle \sigma\right\rangle ^{a}\left\langle \sigma_{1}\right\rangle ^{a_{1}}\left\langle \sigma_{2}\right\rangle ^{a_{2}}\left\langle \xi_{2}\right\rangle ^{k}\left\langle \xi\right\rangle ^{l}}\leq & C\int\frac{\hat{v}\hat{v}_{1}\hat{v}_{2}\left\langle \xi_{1}\right\rangle ^{k-l}}{\left\langle \sigma\right\rangle ^{a}\left\langle \sigma_{1}\right\rangle ^{a_{1}}\left\langle \sigma_{2}\right\rangle ^{a_{2}}\left\langle \xi_{2}\right\rangle ^{k}}=:Z.
\end{align*}
Note that, in this region
\[
\left|\xi_{1}\right|\geq2\left|\xi_{2}\right|\Rightarrow\left|\xi\right|=\left|\xi_{1}-\xi_{2}\right|\geq\left|\left|\xi_{1}\right|-\left|\xi_{2}\right|\right|\geq\left|\xi_{2}\right|.
\]
Moreover, since $\left|\xi_{1}\right|\geq2\left|\xi\right|\Rightarrow-\frac{1}{4}\xi_{1}^{2}\leq-\xi_{2}^{2}$,
we have
\[
\frac{3}{4}\xi_{1}^{2}=\xi_{1}^{2}-\frac{1}{4}\xi_{1}^{2}\leq\xi_{1}^{2}-\xi_{2}^{2}=z\leq\xi_{1}^{2}.
\]
On the other hand, since $\left|\xi_{1}\right|\geq2\left|\xi_{2}\right|\Rightarrow\xi_{1}^{2}\geq4\xi_{2}^{2},$
we have
\[
z=\xi_{1}^{2}-\xi_{2}^{2}\geq3\xi_{2}^{2}.
\]
Summarizing
\begin{equation}
\left|\xi_{1}\right|\geq2\left|\xi_{2}\right|,\textrm{ }\left|\xi\right|\geq\left|\xi_{2}\right|,\textrm{ }\left|\xi_{1}\right|\leq2\left|\xi\right|\leq3\left|\xi_{1}\right|,\label{eq:est_xi_1}
\end{equation}
\begin{equation}
\frac{3}{4}\xi_{1}^{2}\leq z\leq\xi_{1}^{2},\textrm{ }z\geq3\xi_{2}^{2}.\label{eq:est_xi_2}
\end{equation}
Furthermore, it follows from (\ref{eq:Def_z}) and from $\xi=\xi_{1}-\xi_{2}$
that
\begin{align}
z+\xi^{2}= & \xi_{1}^{2}-\xi_{2}^{2}+\xi^{2}=\xi_{1}^{2}-\xi_{2}^{2}+\xi_{1}^{2}-2\xi_{1}\xi_{2}+\xi_{2}^{2}\label{eq:z_xi_xi1}\\
= & 2\xi_{1}^{2}-2\xi_{1}\xi_{2}=2\xi_{1}\left(\xi_{1}-\xi_{2}\right)=2\xi_{1}\xi.\nonumber 
\end{align}
and
\begin{align}
z-\xi^{2}= & \xi_{1}^{2}-\xi_{2}^{2}-\xi^{2}=\xi_{1}^{2}-\xi_{2}^{2}-\xi_{1}^{2}+2\xi_{1}\xi_{2}-\xi_{2}^{2}\nonumber \\
= & 2\xi_{1}\xi_{2}-\xi_{2}^{2}=2\xi_{2}\left(\xi_{1}-\xi_{2}\right)=2\xi_{2}\xi.\label{eq:z_xi_xi2}
\end{align}
And therefore, by (\ref{eq:est_xi_1})
\[
z+\xi^{2}=2\xi\xi_{1}\leq2\left|\xi\right|\left|\xi_{1}\right|\leq4\xi^{2}\Rightarrow z\leq3\xi^{2}.
\]
Moreover, using (\ref{eq:est_xi_1}) and (\ref{eq:est_xi_2}), we
obtain
\[
z+\xi^{2}=2\xi\xi_{1}\leq2\left|\xi\right|\left|\xi_{1}\right|\leq2\left(\frac{3}{2}\left|\xi_{1}\right|\right)\left|\xi_{1}\right|=3\xi_{1}^{2}\leq4z\Rightarrow\frac{1}{3}\xi^{2}\leq z.
\]
Hence,
\begin{equation}
\frac{1}{3}\xi^{2}\leq z\leq3\xi^{2}.\label{eq:bound_z_xi}
\end{equation}
We now estimate $Z$ by the Schwarz method.
\[
Estimates\textrm{ }for\textrm{ }Z
\]

$Region\textrm{ }\sigma_{1}\textrm{ }dominant.$ By exactly the same
computation as in the special case, we obtain in the same way as in
(\ref{eq:C1_sp1}) and (\ref{eq:C1_sp2})
\begin{align*}
C_{1}^{2}\leq & C\underset{\xi_{1},\sigma_{1}}{\sup}\left\langle \xi_{1}\right\rangle ^{2k-2l-4a_{1}}\int_{0}^{\xi_{1}^{2}/4}y^{-1/2}\left\langle y\right\rangle ^{-k}dy\int\left\langle \sigma\right\rangle ^{-2a}\left\langle \sigma_{2}\right\rangle ^{-2a_{2}}d\sigma_{2}\\
\leq & C\underset{\xi_{1}}{\sup}\left\langle \xi_{1}\right\rangle ^{2k-2l-4a_{1}}\int_{0}^{\xi_{1}^{2}/4}y^{-1/2}\left\langle y\right\rangle ^{-k}\left\langle y\right\rangle ^{-\alpha_{1}}dy<\infty,
\end{align*}
provided $k-l\leq2a_{1}$ and $\alpha_{1}>1/2$. The additional factor
$\left\langle y\right\rangle ^{-k}$ in the integral does not provide
any improvement since we need already $\alpha_{\text{1}}>1/2$ in
the special case. The last integral again converges at infinity for
all $k\geq0$ but does not yield any decay in $\xi_{1}$. The condition
$k-l\leq2a_{1}$ corresponds to (\ref{eq:Lem1_c}).

$Region\textrm{ }\sigma_{2}\textrm{ }dominant.$ We use the analog
of (\ref{eq:Def_C}) with fixed $\xi_{2}$ and obtain
\[
C_{2}^{2}=\underset{\xi_{2},\sigma_{2}}{\sup}\left\langle \sigma_{2}\right\rangle ^{-2a_{2}}\left\langle \xi_{2}\right\rangle ^{-2k}\int_{2}\left\langle \xi_{1}\right\rangle ^{2k-2l}\left\langle \sigma\right\rangle ^{-2a}\left\langle \sigma_{1}\right\rangle ^{-2a_{1}}d\xi_{1}d\sigma_{1}.
\]
For fixed $\xi_{2}$, it follows from (\ref{eq:Def_z}) that $dz=2\left|\xi_{1}\right|d\xi_{1}$.
Using (\ref{eq:est_xi_2}) and the fact that $\left|z\right|\leq3\left|\sigma_{2}\right|$
for dominant $\sigma_{2}$ and integrating over $\sigma_{1}$ by the
use of Lemma \ref{lem:bracket_lemma}, we get
\[
C_{2}^{2}\leq C\underset{\xi_{2},\sigma_{2}}{\sup}\left\langle \sigma_{2}\right\rangle ^{-2a_{2}}\left\langle \xi_{2}\right\rangle ^{-2k}\int_{3\xi_{2}^{2}}^{3\left|\sigma_{2}\right|}\left|z\right|^{-1/2}\left\langle z\right\rangle ^{k-l}\left\langle z+\sigma_{2}\right\rangle ^{-\alpha_{2}}dz.
\]
We assume without loss of generality that $k\geq l$. We estimate
the last integral by separating the region $0\leq z\leq\left|\sigma_{2}\right|/2$
and $\left|\sigma_{2}\right|/2\leq z\leq3\left|\sigma_{2}\right|,$
which in the worst case $\sigma_{2}<0$ contribute respectively 
\[
\left\langle \sigma_{2}\right\rangle ^{1/2+k-l-\alpha_{2}},
\]
\[
\left\langle \sigma_{2}\right\rangle ^{-1/2+k-l+\left[1-\alpha_{2}\right]_{+}}.
\]
Keeping the largest contribution, namely the second one, we obtain
\[
C_{2}^{2}\leq C\underset{\sigma_{2}}{\sup}\left\langle \sigma_{2}\right\rangle ^{-2a_{2}-1/2+k-l+\left[1-\alpha_{2}\right]_{+}},
\]
and the last quantity is finite provided
\begin{equation}
k-l\leq2a_{2}+1/2-\left[1-\alpha_{2}\right]_{+}.\label{eq:cond_sig2_dom}
\end{equation}
We shall analyze that condition below together with a similar condition
coming from the region $\sigma$ dominant.

$Region\textrm{ }\sigma\textrm{ }dominant$. We use (\ref{eq:Def_C})
to get
\begin{equation}
C_{*}^{2}=\underset{\xi,\sigma}{\sup}\left\langle \sigma\right\rangle ^{-2a}\left\langle \xi\right\rangle ^{2k-2l}\int_{*}\left\langle \xi_{2}\right\rangle ^{-2k}\left\langle \sigma_{1}\right\rangle ^{-2a_{1}}\left\langle \sigma_{2}\right\rangle ^{-2a_{2}}d\xi_{2}d\sigma_{2}.\label{eq:sig_dom_1}
\end{equation}
Now $\sigma$ dominant implies $\left|z\right|\leq3\left|\sigma\right|$
and therefore $\xi^{2}\leq9\left|\sigma\right|$ by (\ref{eq:bound_z_xi}).
We use this fact to estimate the first factor $\left\langle \sigma\right\rangle ^{-2a}$
in (\ref{eq:sig_dom_1}). It follows again from (\ref{eq:Def_z})
that $dz=2\left|\xi\right|d\xi_{2}$ for fixed $\xi$. We furthermore
express $\xi_{2}$ in terms of $z$ and $\xi$ by (\ref{eq:z_xi_xi2}),
and we integrate over $\sigma_{2}$ for fixed $z$ using Lemma \ref{lem:bracket_lemma}.
We obtain
\[
C_{*}^{2}=C\underset{\xi,\sigma}{\sup}\left\langle \xi\right\rangle ^{2k-2l-4a}\left|\xi\right|^{-1}\int_{\xi^{2}/3}^{3\xi^{2}}\left\langle \left(z-\xi^{2}\right)/2\left|\xi\right|\right\rangle ^{-2k}\left\langle z+\sigma\right\rangle ^{-\alpha}dz.
\]
We next extend the range of integration of $z$ symmetrically to $-2\xi^{2}\leq z-\xi^{2}=y\leq2\xi^{2}$
and apply Lemma \ref{lem:Conv} with $f\left(y\right)=\left\langle y/2\left|\xi\right|\right\rangle ^{-2k}\chi\left(\left|y\right|\leq2\xi^{2}\right)$,
$g\left(y\right)=\left\langle y\right\rangle ^{-\alpha}$ to conclude
that the supremum over $\sigma$ occurs for $\sigma=-\xi^{2}$, so
that
\begin{equation}
C_{*}^{2}=C\underset{\xi}{\sup}\left\langle \xi\right\rangle ^{2k-2l-4a}\left|\xi\right|^{-1}\int_{0}^{2\xi^{2}}\left\langle y/2\left|\xi\right|\right\rangle ^{-2k}\left\langle y\right\rangle ^{-\alpha}dy.\label{eq:estC1}
\end{equation}
The right-hand side of the last expression is bounded for $\left|\xi\right|\leq1$,
i.e., we do not need the restriction $\left|\xi\right|\geq1$. For
$\left|\xi\right|\geq1$ we consider separately the two integration
subregions $0\leq y\leq\left|\xi\right|$ and $\left|\xi\right|\leq y\leq2\xi^{2}$.
The contributions of those regions are estimated respectively by
\begin{equation}
\int_{0}^{\left|\xi\right|}\cdots dy\leq\int_{0}^{\left|\xi\right|}\left\langle y\right\rangle ^{-\alpha}\leq C\left|\xi\right|^{\left[1-\alpha\right]_{+}}dy,\label{eq:estC2}
\end{equation}
\begin{equation}
\int_{\left|\xi\right|}^{2\xi^{2}}\cdots dy\leq C\left|\xi\right|^{2k}\int_{\left|\xi\right|}^{2\xi^{2}}y^{-\alpha-2k}\leq C\left|\xi\right|^{1-\alpha+\left[1-\alpha-2k\right]_{+}}dy.\label{eq:estC3}
\end{equation}
Comparing (\ref{eq:estC1}), (\ref{eq:estC2}), and (\ref{eq:estC3}),
we see that $C_{*}$ is finite provided
\begin{equation}
k-l\leq2a+1/2-\left(1/2\right)\left[1-\alpha\right]_{+},\label{eq:con1_sig_dom}
\end{equation}
\begin{equation}
l>-\left(2a+\alpha\right)+1/2.\label{eq:cond2_sig_dom}
\end{equation}
The last condition holds for any $l\geq-1/2$ provided $2a+\alpha>1$,
which is implied by
\[
a+a_{1}>1/2\textrm{ }a+a_{2}>1/2,\textrm{ }a+a_{1}+a_{2}>1.
\]
Note that the latter set of conditions has already been enforced.
It only remains to ensure (\ref{eq:cond_sig2_dom}) and (\ref{eq:con1_sig_dom}).
Now we have already imposed the conditions $k-l\leq2a_{1}$ and $\alpha_{1}>1/2$,
$\alpha_{2}>1/2$ or equivalently 
\begin{equation}
a,a_{1},a_{2}>1/4,\textrm{ }a+a_{1}>3/4,\textrm{ }a+a_{2}>3/4.\label{eq:cond_a_alphas}
\end{equation}
The conditions (\ref{eq:cond_sig2_dom}) and (\ref{eq:con1_sig_dom})
are implied respectively by
\begin{align}
k-l< & 2a_{2}+1/2,\label{eq:c1_a}\\
k-l< & 2a_{2}+2a-1/2,\label{eq:c1_b}\\
k-l< & 2a_{2}+2a_{1}-1/2,\label{eq:c1_c}\\
k-l< & 2a_{2}+2a+2a_{1}-3/2,\label{eq:c1_d}
\end{align}
and
\begin{align}
k-l< & 2a+1/2,\label{eq:c2_a}\\
k-l< & 2a+a_{1},\label{eq:c2_b}\\
k-l< & 2a+a_{2},\label{eq:c2_c}\\
k-l< & 2a+a_{1}+a_{2}-1/2.\label{eq:c2_d}
\end{align}
Now $k-l\leq2a_{1}$ and (\ref{eq:cond_a_alphas}) imply (\ref{eq:c1_c})
and (\ref{eq:c1_d}). Next, $2a+a_{1}>a+1/4+a_{1}=\left(1/2\right)\left(2a+1/2+2a_{1}\right)$
so that $k-l\leq2a_{1}$ and (\ref{eq:c2_a}) imply (\ref{eq:c2_b}).
Furthermore, $2a+a_{2}=\left(1/2\right)\left(2a_{2}+2a-1/2+2a+1/2\right)$
so that (\ref{eq:c1_b}) and (\ref{eq:c2_a}) imply (\ref{eq:c2_c}).
Finally, $2a+a_{1}+a_{2}-1/2>a+a_{1}+a_{2}-1/4=\left(1/2\right)\left(2a_{1}+2a_{2}+2a-1/2\right)$
so that $k-l\leq2a_{1}$ and (\ref{eq:c1_b}) imply (\ref{eq:c2_d}).
It is therefore sufficient to ensure (\ref{eq:c1_a}), (\ref{eq:c1_b}),
and (\ref{eq:c2_a}). By (\ref{eq:cond_a_alphas}), the right-hand
side of those three inequalities are all $>1$; they are implied by
$k-l\leq1$, contained in (\ref{eq:Lem1_a}).
\end{proof}
Now we have to verify the bounds for $A$ in (\ref{eq:s_leq_A_Z0}),
that is
\begin{equation}
\int\frac{\left|\hat{v}\hat{v}_{1}\hat{v}_{2}\right|}{\left\langle \sigma\right\rangle ^{a}\left\langle \sigma_{1}\right\rangle ^{a_{1}}\left\langle \sigma_{2}\right\rangle ^{a_{2}}}\leq C\left\Vert v\right\Vert _{2}\left\Vert v_{1}\right\Vert _{2}\left\Vert v_{2}\right\Vert _{2},\label{eq:ver_A}
\end{equation}
provided $a,a_{1},a_{2}>1/4$ (see the first condition in (\ref{eq:Lem1_b})).
We use the following result for the Schrödinger equation (see \cite[Lemma 2.4]{ginibre1997cauchy})
\begin{lem}
\label{lem:Lemm_Sch}Let $\phi\left(\xi\right)=\xi^{2}$ (Schrödinger
equation). Assume $b_{0}>1/2$ , $0\leq b\leq b_{0}$, and $0<\eta\leq1$
($\eta\geq1/2$ if $n=1$, i.e., 1D). Then
\[
\left\Vert f\right\Vert _{L_{t}^{q}\left(L_{x}^{r}\right)}\leq C\left\Vert f\right\Vert _{X_{\phi\left(\xi\right)=\xi^{2}}^{0,b}},
\]
where $2/q=1-\eta b/b_{0}$, $\delta\left(r\right)\equiv n/2-n/r=\left(1-\eta\right)b/b_{0}$.
\end{lem}

Using the previous lemma, we show the following
\begin{lem}
\label{lem:Lemm_Sch2}Let $a>1/4,$ $b_{0}=2a>1/2$ . Consider $n=1$.
Let $v\in L^{2}$ and $\alpha=\tau+\xi^{2}$ (Schrödinger). Then
\[
\left\Vert \mathcal{F}^{-1}\left(\left\langle \alpha\right\rangle ^{-a}\left|\hat{v}\right|\right)\right\Vert _{L_{t}^{8/3}\left(L_{x}^{4}\right)}\leq C\left\Vert v\right\Vert _{2}.
\]
\end{lem}

\begin{proof}
By Lemma \ref{lem:Lemm_Sch} with $b=a\leq b_{0}=2a$, $\hat{f}=\left\langle \alpha\right\rangle ^{-a}\left|\hat{v}\right|,$
and $\eta=1/2,$ we have
\[
\left\Vert \mathcal{F}^{-1}\left(\left\langle \alpha\right\rangle ^{-a}\left|\hat{v}\right|\right)\right\Vert _{L_{t}^{q}\left(L_{x}^{r}\right)}\leq C\left\Vert \hat{v}\right\Vert _{2},
\]
where $q=8/3$ and $r=4$,
\end{proof}
\begin{lem}
\label{lem:Lemm_ver_A}Let $a,a_{1},a_{2}>1/4$ and $v,v_{1},v_{2}\in L^{2}$.
Then (\ref{eq:ver_A}) holds.
\end{lem}

\begin{proof}
Since (\ref{eq:ver_A}) is decreasing in $a,a_{1},a_{2},$ it is sufficient
to consider $a=a_{1}=a_{2}>1/4$. We apply Hölder's inequality in
space and time to obtain
\begin{align}
\int\frac{\left|\hat{v}\hat{v}_{1}\hat{v}_{2}\right|}{\left\langle \sigma\right\rangle ^{a}\left\langle \sigma_{1}\right\rangle ^{a_{1}}\left\langle \sigma_{2}\right\rangle ^{a_{2}}}\leq & \left\Vert \mathcal{F}^{-1}\left(\left\langle \sigma\right\rangle ^{-a}\left|\hat{v}\right|\right)\right\Vert _{L_{t}^{q}\left(L_{x}^{r}\right)}\times\nonumber \\
 & \times\prod_{i=1,2}\left\Vert \mathcal{F}^{-1}\left(\left\langle \sigma_{i}\right\rangle ^{-a_{i}}\left|\hat{v}_{i}\right|\right)\right\Vert _{L_{t}^{q_{i}}\left(L_{x}^{r_{i}}\right),}\label{eq:norms_verA}
\end{align}
with
\[
\frac{1}{q}+\frac{1}{q_{1}}+\frac{1}{q_{2}}=1,\textrm{ }\frac{1}{r}+\frac{1}{r_{1}}+\frac{1}{r_{2}}=1.
\]
Let
\[
\left(q,r\right)=\left(4,2\right),\textrm{ }\left(q_{i},r_{i}\right)=\left(\frac{8}{3},4\right),\textrm{ }i=1,2.
\]
Then, the last two norms in (\ref{eq:norms_verA}) are estimated by
Lemma \ref{lem:Lemm_Sch2}. Now recalling the definition of $\sigma$,
we use the Hardy-Littlewood-Sobolev inequality in time to get
\[
\left\Vert \mathcal{F}^{-1}\left(\left\langle \sigma\right\rangle ^{-a}\left|\hat{v}\right|\right)\right\Vert _{L_{t}^{q}\left(L_{x}^{r}\right)}\leq C\left\Vert v\right\Vert _{2},
\]
since $r=2$ and 
\[
\frac{1}{2}-\frac{1}{q}=a,
\]
with $q=4,$ $a=1/4$.
\end{proof}
Our next step is to estimate the nonlinearity $f_{1}=nu$ in $Y_{\phi\left(\xi\right)=\xi^{2}}^{k}$.
For this, we divide $\left|\hat{f}_{1}\right|$ by $\left\langle \sigma_{1}\right\rangle $,
integrate over $\tau_{1}$ (or $\sigma_{1}$) for fixed $\xi_{1}$
and then take the scalar product with a generic function in $H_{x}^{-k}$
with Fourier transform $\left\langle \xi_{1}\right\rangle ^{k}\hat{w}_{1},$
$w_{1}\in L_{x}^{2}$. The estimate of $f_{1}$ in $Y_{\phi\left(\xi\right)=\xi^{2}}^{k}$
becomes
\begin{equation}
\tilde{S}\leq C\left\Vert v\right\Vert _{2}\left\Vert w_{1}\right\Vert _{2}\left\Vert v_{2}\right\Vert _{2},\label{eq:est_S_tilde}
\end{equation}
where
\[
\tilde{S}=\int\frac{\left|\hat{v}\hat{w}_{1}\hat{v}_{2}\right|\left\langle \xi_{1}\right\rangle ^{k}}{\left\langle \sigma\right\rangle ^{a}\left\langle \sigma_{1}\right\rangle \left\langle \sigma_{2}\right\rangle ^{a_{2}}\left\langle \xi_{2}\right\rangle ^{k}\left\langle \xi\right\rangle ^{l}},
\]
with the same notation as in (\ref{eq:def_S}). Note that $w_{1}$
is a function of space only, whereas the $v'$s are functions of space
and time.
\begin{lem}
\label{lem:Lemm_S_tilde}Let $a,$ $a_{2}$, $k$, and $l$ satisfy
(\ref{eq:Lem1_a}) and
\[
a,a_{2}>1/4,\textrm{ }a+a_{2}>3/4.
\]
Then (\ref{eq:est_S_tilde}) holds.
\end{lem}

\begin{proof}
The proof is similar to the one of Lemma \ref{lem:Lemm_S}. However,
we have to handle $w_{1}$ appropriately. For this, let $a_{1}$ satisfy
\[
0<1/2-a_{1}<\min\left(1/4,a-1/4,a+a_{2}-3/4\right),
\]
so that
\[
a_{1}>1/4,\textrm{ }a_{1}+a>3/4,\textrm{ }a_{1}+a_{2}>3/4,\textrm{ }a+a_{1}+a_{2}>5/4.
\]
Define $\hat{v}_{1}=\left\langle \sigma_{1}\right\rangle ^{a_{1}-1}\hat{w}_{1}$.
It follows that $\left\Vert v_{1}\right\Vert _{2}\leq C\left(1-2a_{1}\right)^{-1/2}\left\Vert w_{1}\right\Vert _{2}$.
Under these conditions, one can follow the proof of Lemma \ref{lem:Lemm_S}
with just minor modifications (cf. \cite[Lemma 4.5]{ginibre1997cauchy}).
\end{proof}
Now we consider the nonlinearity $f_{2}=\left|u\right|^{2}u$, which
has been extensibly studied in the context of the NLS equation. We
want to estimate $f_{2}$ in $X_{\phi\left(\xi\right)=\xi^{2}}^{k,-a_{1}}$
. Hence, we have to verify the expression
\begin{equation}
\left|S_{0}\right|\leq C\prod_{i=1}^{4}\left\Vert v_{i}\right\Vert _{2},\label{eq:cub_est1}
\end{equation}
with
\begin{equation}
S_{0}=\int\frac{\hat{v}_{1}\hat{v}_{2}\hat{v}_{3}\hat{v}_{4}\left\langle \xi_{1}\right\rangle ^{k}}{\left\langle \sigma_{1}\right\rangle ^{a_{1}}\left\langle \sigma_{2}\right\rangle ^{a_{2}}\left\langle \sigma_{3}\right\rangle ^{a_{2}}\left\langle \sigma_{4}\right\rangle ^{a_{2}}\left\langle \xi_{2}\right\rangle ^{k}\left\langle \xi_{3}\right\rangle ^{k}\left\langle \xi_{4}\right\rangle ^{k}},\label{eq:cub_est1b}
\end{equation}
where $\hat{v}_{i}=\hat{v}_{i}\left(\xi_{i},\tau_{i}\right),$ $\sigma_{i}=\tau_{i}+\xi_{i}^{2},$
$1\leq i\leq4$. The integral is over $\left(\xi_{i},\tau_{i}\right)$,
constrained by $\xi_{1}+\xi_{2}=\xi_{3}+\xi_{4}$ and $\tau_{1}+\tau_{2}=\tau_{3}+\tau_{4}$.
Furthermore, if either $a_{1}=1/2$ or $a_{2}\leq1/2$, we need to
estimate $f_{2}$ in $Y_{\phi\left(\xi\right)=\xi^{2}}^{k}$, hence,
we have to verify 
\begin{equation}
\left|\tilde{S}_{0}\right|\leq C\left\Vert w_{1}\right\Vert _{2}\prod_{i=2}^{4}\left\Vert v_{i}\right\Vert _{2},\label{eq:cub_est2}
\end{equation}
with
\begin{equation}
\tilde{S}_{0}=\int\frac{\hat{w}_{1}\hat{v}_{2}\hat{v}_{3}\hat{v}_{4}\left\langle \xi_{1}\right\rangle ^{k}}{\left\langle \sigma_{1}\right\rangle \left\langle \sigma_{2}\right\rangle ^{a_{2}}\left\langle \sigma_{3}\right\rangle ^{a_{2}}\left\langle \sigma_{4}\right\rangle ^{a_{2}}\left\langle \xi_{2}\right\rangle ^{k}\left\langle \xi_{3}\right\rangle ^{k}\left\langle \xi_{4}\right\rangle ^{k}},\label{eq:cub_est2b}
\end{equation}
where we are using the same notation as before.
\begin{lem}
\label{lem:Lemm_cub}Let $k\geq0$ and 
\[
\max\left(1/6,\left(1-k\right)/3\right)<a_{2}<1.
\]
Then, (\ref{eq:cub_est1}) and (\ref{eq:cub_est2}) hold.
\end{lem}

See \cite[Lemma 4.7]{ginibre1997cauchy} for a sketch of the proof.
See, e.g., \cite{bourgain1993fourier,kenig1996quadratic} for additional
details. 

Next we consider the nonlinearity $f_{3}=\left|u\right|^{2}n$. We
want to estimate $f_{3}$ in $X_{\phi\left(\xi\right)\equiv0}^{l,-a_{0}}$,
for suitable $a_{0}$. Hence, we have to verify the expression
\begin{equation}
\left|S_{1}\right|\leq C\prod_{i=1}^{4}\left\Vert v_{i}\right\Vert _{2},\label{eq:uun_est}
\end{equation}
with
\begin{equation}
S_{1}=\int\frac{\hat{v}_{1}\hat{v}_{2}\hat{v}_{3}\hat{v}_{4}\left\langle \xi_{1}\right\rangle ^{l}}{\left\langle \sigma_{1}\right\rangle ^{a_{0}}\left\langle \sigma_{2}\right\rangle ^{a_{2}}\left\langle \sigma_{3}\right\rangle ^{a_{2}}\left\langle \sigma_{4}\right\rangle ^{a}\left\langle \xi_{2}\right\rangle ^{k}\left\langle \xi_{3}\right\rangle ^{k}\left\langle \xi_{4}\right\rangle ^{l}},\label{eq:uun_est1b}
\end{equation}
where $\hat{v}_{i}=\hat{v}_{i}\left(\xi_{i},\tau_{i}\right)\textrm{ }1\leq i\leq4,$
$\sigma_{1}=\tau_{1}$, $\sigma_{2}=\tau_{2}+\xi_{2}^{2},$ $\sigma_{3}=\tau_{3}+\xi_{3}^{2}$,
$\sigma_{4}=\tau_{4}$. The integral is over $\left(\xi_{i},\tau_{i}\right)$,
constrained by $\xi_{1}+\xi_{2}=\xi_{3}+\xi_{4}$ and $\tau_{1}+\tau_{2}=\tau_{3}+\tau_{4}$.
Note that $\hat{v}_{2}\left(\xi_{2},\tau_{2}\right)$ implies that
$\hat{\bar{v}}_{2}\left(-\xi_{2},-\tau_{2}\right)$, the change of
sign due to complex conjugation. Furthermore, if either $a_{0}=1/2$
or $a\leq1/2$, we need to estimate $f_{3}$ in $Y_{\phi\left(\xi\right)\equiv0}^{k}$,
hence, we have to verify 
\begin{equation}
\left|\tilde{S}_{1}\right|\leq C\left\Vert w_{1}\right\Vert _{2}\prod_{i=2}^{4}\left\Vert v_{i}\right\Vert _{2},\label{eq:uun_est2}
\end{equation}
with
\begin{equation}
\tilde{S}_{1}=\int\frac{\hat{w}_{1}\hat{v}_{2}\hat{v}_{3}\hat{v}_{4}\left\langle \xi_{1}\right\rangle ^{l}}{\left\langle \sigma_{1}\right\rangle \left\langle \sigma_{2}\right\rangle ^{a_{2}}\left\langle \sigma_{3}\right\rangle ^{a_{2}}\left\langle \sigma_{4}\right\rangle ^{a}\left\langle \xi_{2}\right\rangle ^{k}\left\langle \xi_{3}\right\rangle ^{k}\left\langle \xi_{4}\right\rangle ^{l}},\label{eq:uun_est2b}
\end{equation}
where we are using the same notation as before. We need the following
intermediate result.
\begin{lem}
\label{lem:Int_u2n_est}Let 
\[
W:=\int\frac{\hat{v}_{1}\hat{v}_{2}\hat{v}_{3}\hat{v}_{4}}{\left\langle \sigma_{1}\right\rangle ^{a_{0}}\left\langle \sigma_{2}\right\rangle ^{a_{2}}\left\langle \sigma_{3}\right\rangle ^{a_{2}}\left\langle \sigma_{4}\right\rangle ^{a}\left\langle \xi_{2}\right\rangle ^{k}\left\langle \xi_{3}\right\rangle ^{k}\left\langle \xi_{4}\right\rangle ^{l}}.
\]
Then
\[
\left|W\right|\leq C\prod_{i=1}^{4}\left\Vert v_{i}\right\Vert _{2},
\]
provided
\begin{equation}
l\geq\delta_{4}=1-2\left(\left(1-\eta\right)\frac{a_{2}}{b_{0}}+k\right),\textrm{ }l>\frac{1}{2}\textrm{ if }\delta_{4}=\frac{1}{2},\label{eq:LL1}
\end{equation}
\begin{equation}
a_{0}+a+\eta\frac{a_{2}}{b_{0}}=1,\label{eq:LL2}
\end{equation}
with $1/2\leq\eta\leq1$, $b_{0}\geq a_{2},$ $b_{0}\geq1/2$. In
particular, for $\eta=1/2$ and $b_{0}=a_{2}>1/2$, we require
\begin{equation}
l\geq-2k,\textrm{ }a_{0}+a=1/2,\textrm{ and }a_{2}>1/2.\label{eq:LL3}
\end{equation}
\end{lem}

\begin{proof}
Using Hölder's inequality in space and time, we have
\begin{align*}
\left|W\right|\leq & \left\Vert \mathcal{F}^{-1}\left(\left\langle \sigma_{1}\right\rangle ^{-a_{0}}\left|\hat{v}_{1}\right|\right)\right\Vert _{L_{t}^{q_{1}}\left(L_{x}^{r_{1}}\right)}\times\\
 & \times\prod_{i=2,3}\left\Vert \mathcal{F}^{-1}\left(\left\langle \xi_{i}\right\rangle ^{-k}\left\langle \sigma_{i}\right\rangle ^{-a_{2}}\left|\hat{v}_{i}\right|\right)\right\Vert _{L_{t}^{q_{i}}\left(L_{x}^{r_{i}}\right)}\times\\
 & \times\left\Vert \mathcal{F}^{-1}\left(\left\langle \xi_{4}\right\rangle ^{-l}\left\langle \sigma_{4}\right\rangle ^{-a}\left|\hat{v}_{4}\right|\right)\right\Vert _{L_{t}^{q_{4}}\left(L_{x}^{r_{4}}\right)},
\end{align*}
with
\begin{equation}
\frac{1}{q_{1}}+\frac{1}{q_{2}}+\frac{1}{q_{3}}+\frac{1}{q_{4}}=1,\label{eq:Holder1}
\end{equation}

\begin{equation}
\frac{1}{r_{1}}+\frac{1}{r_{2}}+\frac{1}{r_{3}}+\frac{1}{r_{4}}=1,\textrm{ or }\delta_{1}+\delta_{2}+\delta_{3}+\delta_{4}=1,\label{eq:Holder2}
\end{equation}
and $r_{1}=2$. Using an argument similar to Lemma \ref{lem:Lemm_Sch2}
and the Hardy-Littlewood-Sobolev inequality, we obtain
\[
\left|W\right|\leq C\prod_{i=1}^{4}\left\Vert v_{i}\right\Vert _{2},
\]
provided
\begin{equation}
\frac{2}{q_{1}}=1-2a_{0},\textrm{ }r_{1}=2\Rightarrow\delta_{1}=0,\label{eq:cc1}
\end{equation}
\begin{equation}
\frac{2}{q_{2}}=1-\eta\frac{a_{2}}{b_{0}},\textrm{ }H^{k,\tilde{r}_{2}}\subset L^{r_{2}},\textrm{ }\tilde{\delta}_{2}=\left(1-\eta\right)\frac{a_{2}}{b_{0}},\label{eq:cc2}
\end{equation}
\begin{equation}
\frac{2}{q_{3}}=1-\eta\frac{a_{2}}{b_{0}},\textrm{ }H^{k,\tilde{r}_{3}}\subset L^{r_{3}},\textrm{ }\tilde{\delta}_{3}=\left(1-\eta\right)\frac{a_{2}}{b_{0}},\label{eq:cc3}
\end{equation}
\begin{equation}
\frac{2}{q_{4}}=1-2a,\textrm{ }H^{l}\subset L^{r_{4}},\textrm{ }l\geq\delta_{4}=\frac{1}{2}-\frac{1}{r_{4}}\geq0.\label{eq:cc4}
\end{equation}
and $\frac{1}{2}\leq\eta\leq1$. 

Considering (\ref{eq:Holder1}) along with the $q's$ in (\ref{eq:cc1})-(\ref{eq:cc4}),
we get (\ref{eq:LL2}). From (\ref{eq:cc1}), we have
\[
\frac{1}{r_{2}}=\frac{1}{\tilde{r}_{2}}-k,\textrm{ }k<\frac{1}{\tilde{r}_{2}}\Rightarrow\tilde{\delta}_{2}=\frac{1}{2}-\frac{1}{r_{2}}-k=\delta_{2}-k.
\]
Hence,
\[
\delta_{2}=\left(1-\eta\right)\frac{a_{2}}{b_{0}}+k.
\]
Similarly,
\[
\delta_{3}=\left(1-\eta\right)\frac{a_{2}}{b_{0}}+k.
\]
Combining the expressions for the $\delta's$ with the RHS of (\ref{eq:Holder2}),
we obtain (\ref{eq:LL1}). (\ref{eq:LL3}) follows from the previous
results.
\end{proof}
\begin{lem}
\label{lem:lemm_cub_ode}Let 
\[
k\geq0,\textrm{ }l\leq k,\textrm{ }a_{0}+a=1/2,\textrm{ }a_{2}>1/2.
\]
Then (\ref{eq:uun_est}) and (\ref{eq:uun_est2}) hold.
\end{lem}

\begin{proof}
Note that, due to the constraint $\xi_{1}+\xi_{2}=\xi_{3}+\xi_{4}$,
we get $\left\langle \xi_{1}\right\rangle ^{l}\leq C\left(\left\langle \xi_{2}\right\rangle ^{l}+\left\langle \xi_{3}\right\rangle ^{l}+\left\langle \xi_{4}\right\rangle ^{l}\right)$.
Then, considering the symmetry in the variables 2 and 3, we obtain
\begin{align*}
\left|S_{1}\right|\leq & C\int\frac{\left|\hat{v}_{1}\hat{v}_{2}\hat{v}_{3}\hat{v}_{4}\right|\left\langle \xi_{2}\right\rangle ^{l-k}}{\left\langle \sigma_{1}\right\rangle ^{a_{0}}\left\langle \sigma_{2}\right\rangle ^{a_{2}}\left\langle \sigma_{3}\right\rangle ^{a_{2}}\left\langle \sigma_{4}\right\rangle ^{a}\left\langle \xi_{3}\right\rangle ^{k}\left\langle \xi_{4}\right\rangle ^{l}}+\\
 & +C\int\frac{\left|\hat{v}_{1}\hat{v}_{2}\hat{v}_{3}\hat{v}_{4}\right|}{\left\langle \sigma_{1}\right\rangle ^{a_{0}}\left\langle \sigma_{2}\right\rangle ^{a_{2}}\left\langle \sigma_{3}\right\rangle ^{a_{2}}\left\langle \sigma_{4}\right\rangle ^{a}\left\langle \xi_{2}\right\rangle ^{k}\left\langle \xi_{3}\right\rangle ^{k}}:=A+B.
\end{align*}
To bound $B$, we use Lemma \ref{lem:Int_u2n_est} with $l=0$. Hence,
we require
\[
0\geq-2k,\textrm{ }a_{0}+a=1/2,\textrm{ }a_{2}>1/2,
\]
which implies $k\geq0$. For $A$, we use the condition $l-k\leq0$.
Then, we consider Lemma \ref{lem:Int_u2n_est} replacing $k\mapsto k-l$,
since $A$ is increasing in $l$ and decreasing in $k$ (recall $k\geq0$,
$k\geq l$). Then, we need
\[
-l\geq-2k,\textrm{ }a_{0}+a=1/2,\textrm{ }a_{2}>1/2.
\]
Hence, combining all the previous conditions we obtain the result
for (\ref{eq:uun_est}). 

To estimate $\tilde{S}_{1}$, take $\hat{v}_{1}=\left\langle \sigma_{1}\right\rangle ^{a_{0}-1}\hat{w}_{1}$
in (\ref{eq:uun_est1b}). Since for $a_{0}<1/2$, we have $\left\Vert \hat{v}_{1}\right\Vert _{2}\leq C\left(a_{0}\right)\left\Vert \hat{w}_{1}\right\Vert _{2}$,
(\ref{eq:uun_est2}) holds.
\end{proof}
We shall use the following simple observation.
\begin{lem}
\label{lem:lemm_aux}For any $\varepsilon>0$
\begin{equation}
\left\Vert u\right\Vert _{Y_{\phi\left(\xi\right)}^{s}}\leq C\left\Vert u\right\Vert _{X_{\phi\left(\xi\right)}^{s,-1/2+\varepsilon}.}\label{eq:Ys_Xsb}
\end{equation}
\end{lem}

\begin{proof}
Applying the Cauchy-Schwarz inequality in $\tau$, we get
\begin{align*}
\left\Vert u\right\Vert _{Y_{\phi\left(\xi\right)}^{s}}^{2}= & \int\left(\left\langle \xi\right\rangle ^{s}\int\left\langle \tau+\phi\left(\xi\right)\right\rangle ^{-1}\left|\hat{u}\left(\tau,\xi\right)\right|d\tau\right)^{2}d\xi\\
\leq & \int\left\langle \xi\right\rangle ^{2s}\left(\int\left\langle \tau+\phi\left(\xi\right)\right\rangle ^{2\left(-\frac{1}{2}-\varepsilon\right)}d\tau\right)\left(\int\left\langle \tau+\phi\left(\xi\right)\right\rangle ^{2\left(-\frac{1}{2}+\varepsilon\right)}\left|\hat{u}\left(\tau,\xi\right)\right|^{2}d\tau\right)d\xi,
\end{align*}
for any $\varepsilon>0$. Since $-\left(1+2\varepsilon\right)<-1$,
(\ref{eq:Ys_Xsb}) follows. 
\end{proof}
Now we present the main result of this section
\begin{prop}
The exciton-polariton system (\ref{eq:SYS}) with initial data $\left(u_{0},n_{0}\right)\in H^{k}\oplus H^{l}$
is locally well-posed in $X_{\phi\left(\xi\right)=\xi^{2},T}^{k,a_{2}}\oplus X_{\phi\left(\xi\right)\equiv0,T}^{l,a}$
provided
\begin{equation}
k\geq0,\textrm{ }l\leq k,\textrm{ }k-l\leq2a_{1},\label{eq:c1}
\end{equation}
\begin{equation}
a=1/4+3\varepsilon,\label{eq:c2}
\end{equation}
\begin{equation}
a_{1}=1/2-2\varepsilon,\label{eq:c3}
\end{equation}
\begin{equation}
a_{2}=1/2+\varepsilon,\label{eq:c4}
\end{equation}
with $\varepsilon>0$ small enough ($\varepsilon<1/12$). Moreover,
\[
\left(u,n\right)\in C\left(\left[-T,T\right];H^{k}\oplus H^{l}\right),
\]
with $T=T\left(\left\Vert u_{0}\right\Vert _{H^{k}},\left\Vert n_{0}\right\Vert _{H^{l}}\right)>0$.
\end{prop}

\begin{proof}
Set 
\begin{equation}
a_{0}=1/4-3\varepsilon.\label{eq:c5}
\end{equation}
Then, one can verify that under conditions (\ref{eq:c1})-(\ref{eq:c5})
all the assumptions of Lemmas \ref{lem:Lemm_S}, \ref{lem:Lemm_S_tilde},
\ref{lem:Lemm_cub}, \ref{lem:lemm_cub_ode}, \ref{lem:lemm_aux}
are satisfied. Moreover, we have $a+a_{0}<1$ and $a_{2}+a_{1}<1$,
hence, we can apply Lemma \ref{lem:Lem_ini} and obtain therefrom
a strictly positive power of $T$. Then, we get the result considering
the cutoff system (\ref{eq:EP1_cut})-(\ref{eq:EP2_cut}) and using
a standard fixed point argument. Notice that we use the spaces restricted
in time to deal with the term $P=P\left(x\right)$, not to get a positive
power of $T$, which we get from Lemma \ref{lem:Lem_ini}. For $a\leq1/2$
we have to take into account \cite[Lemma 2.2]{ginibre1997cauchy}
to conclude continuity in time of the solution. 
\end{proof}
\begin{cor}
Let $\left(u_{0},n_{0}\right)\in L^{2}\oplus L^{2}$ with $n_{0}\left(x\right)\geq0$.
Then, there exists a global in time solution $\left(u,n\right)\in C\left(\left[0,\infty\right),L^{2}\oplus L^{2}\right)$
of the exciton-polariton system (\ref{eq:SYS}). Furthermore, the
system has an absorbing set in $L^{2}\oplus L^{2}$.
\end{cor}

\begin{proof}
Consider a smooth solution of (\ref{eq:SYS}), then argue by density.
Using the usual variation of constants formula in the second equation
of (\ref{eq:SYS}), we have 
\[
n\left(t,x\right)=n_{0}\left(x\right)e^{-\int_{0}^{t}\Gamma\left(\tau,x\right)d\tau}+P\int_{0}^{t}e^{-\int_{s}^{t}\Gamma\left(\tau,x\right)d\tau}ds,
\]
where $\Gamma\left(t,x\right)=R\left|u\left(t,x\right)\right|^{2}+\beta$.
Hence, if $n_{0}\left(x\right)\geq0$, then $n\left(t,x\right)\geq0$
for all $t\in\left[0,T\right]$ since $P=P\left(x\right)\geq0$. 

Now multiply the first equation in (\ref{eq:SYS}) by $\bar{u}$,
integrate over $\mathbb{R},$ and take the imaginary part to get
\[
\frac{d}{dt}\frac{1}{2}\int_{\mathbb{R}}\left|u\right|^{2}dx=\int_{\mathbb{R}}\left(Rn-\alpha\right)\left|u\right|^{2}dx.
\]
Furthermore, integrate the second equation in (\ref{eq:SYS}) over
$\mathbb{R}$ to obtain
\[
\frac{d}{dt}\int_{\mathbb{R}}ndx=\int\left[P-\left(R\left|u\right|^{2}+\beta\right)n\right]dx.
\]
Combining the last two expressions gives
\begin{align*}
\frac{d}{dt}\left(\frac{1}{2}\int_{\mathbb{R}}\left|u\right|^{2}dx+\int_{\mathbb{R}}ndx\right) & =\int_{\mathbb{R}}\left[\left(Rn-\alpha\right)\left|u\right|^{2}+P-\left(R\left|u\right|^{2}+\beta\right)n\right]dx\\
 & =\int_{\mathbb{R}}Pdx-\alpha\int_{\mathbb{R}}\left|u\right|^{2}dx-\beta\int_{\mathbb{R}}ndx\\
 & \leq\int_{\mathbb{R}}Pdx-\gamma\left(\frac{1}{2}\int_{\mathbb{R}}\left|u\right|^{2}dx+\int_{\mathbb{R}}ndx\right),
\end{align*}
where $\gamma=\min\left(2\alpha,\beta\right)$. Integrating in time,
we get
\begin{equation}
\frac{1}{2}\int_{\mathbb{R}}\left|u\right|^{2}dx+\int_{\mathbb{R}}ndx\leq e^{-\gamma t}\left(\frac{1}{2}\int_{\mathbb{R}}\left|u_{0}\right|^{2}dx+\int_{\mathbb{R}}n_{0}dx-\frac{1}{\gamma}\int_{\mathbb{R}}Pdx\right)+\frac{1}{\gamma}\int_{\mathbb{R}}Pdx.\label{eq:Mass1}
\end{equation}

Now multiply the second equation in (\ref{eq:SYS}) by $2n$ to obtain
\[
\partial_{t}n^{2}=2Pn-2\left(R\left|u\right|^{2}+\beta\right)n^{2}\leq\frac{P^{2}}{\beta}-\beta n^{2},
\]
where the last inequality follows from $\left(\frac{P}{\sqrt{\beta}}-\sqrt{\beta}n\right)^{2}\geq0.$
This implies that
\[
\partial_{t}\left(e^{t\beta}n^{2}\right)\leq e^{t\beta}\frac{P^{2}}{\beta}.
\]
Integrating the last expression in time, we get
\[
n^{2}\left(\cdot,t\right)\leq e^{-t\beta}\left(n_{0}^{2}-\frac{P^{2}}{\beta^{2}}\right)+\frac{P^{2}}{\beta^{2}},\textrm{ }\forall0\leq t\leq T.
\]
Hence,
\begin{equation}
\int_{\mathbb{R}}n^{2}dx\leq e^{-t\beta}\left(\int_{\mathbb{R}}n_{0}dx-\frac{1}{\beta^{2}}\int_{\mathbb{R}}P^{2}dx\right)+\frac{1}{\beta}\int_{\mathbb{R}}P^{2}dx.\label{eq:mass2}
\end{equation}

The result follows by combining (\ref{eq:Mass1}), (\ref{eq:mass2}),
the fact that $n_{0}\left(x\right)\geq0$, and a density argument.
\end{proof}
The global existence theory of (\ref{eq:SYS}) in $H^{1}\left(\mathbb{T}\right)\oplus H^{1}\left(\mathbb{T}\right)$
was established in \cite{antonelli2019dissipative}.

\section{Acknowledgement }

The authors acknowledge financial support from the Austrian Science
Fund (FWF) grant F65.

\section{References }

\bibliographystyle{plain}
\bibliography{bib}

\end{document}